%% file: msimp.tex
\documentclass{amsart}
\input{aux/style}
\input{aux/usualcmds}

\input{aux/commands}

\title[Multisimplicial chains and configuration spaces]{Multisimplicial chains and configuration spaces}

\author[Medina-Mardones]{Anibal M. Medina-Mardones}
\address{A.M-M., LAGA, Universit\'e Sorbonne Paris Nord}
\email{\href{medina-mardones@math.univ-paris13.fr}{medina-mardones@math.univ-paris13.fr}}

\author[Pizzi]{Andrea Pizzi}
\address{A.P., Dipartimento di Matematica, Universit\`a di Roma Tor Vergata}
\email{\href{pizzi@mat.uniroma2.it}{pizzi@mat.uniroma2.it}}

\author[Salvatore]{Paolo Salvatore}
\address{P.S., Dipartimento di Matematica, Universit\`a di Roma Tor Vergata}
\email{\href{salvator@mat.uniroma2.it}{salvator@mat.uniroma2.it}}

\date{\today}

\begin{document}
	\input{sec/abstract}
	\maketitle
	\input{sec/introduction}
	\input{sec/acknowledgment}
	\input{sec/multisimplicial}
	\input{sec/comparison}
	\input{sec/configuration}
	\sloppy
	\printbibliography
\end{document}

%% file: aux/style.tex
\usepackage{microtype}
\usepackage{amssymb}
\usepackage{mathtools}
\usepackage{tikz-cd}
\usepackage{mathbbol} % changes \mathbb{} and adds more support

\usepackage[
	backend=biber,
	style=alphabetic,
	backref=true,
	url=false,
	doi=false,
	isbn=false,
	eprint=false]{biblatex}

\renewbibmacro{in:}{} % don't display "in:" before the journal name
\AtEveryBibitem{\clearfield{pages}} % don't show page numbers

\DeclareFieldFormat{title}{\myhref{\mkbibemph{#1}}}
\DeclareFieldFormat
[article,inbook,incollection,inproceedings,patent,thesis,unpublished]
{title}{\myhref{\mkbibquote{#1\isdot}}}

\newcommand{\doiorurl}{%
	\iffieldundef{url}
	{\iffieldundef{eprint}
		{}
		{http://arxiv.org/abs/\strfield{eprint}}}
	{\strfield{url}}%
}

\newcommand{\myhref}[1]{%
	\ifboolexpr{%
		test {\ifhyperref}
		and
		not test {\iftoggle{bbx:eprint}}
		and
		not test {\iftoggle{bbx:url}}
	}
	{\href{\doiorurl}{#1}}
	{#1}%
}

\usepackage[
	bookmarks=true,
	linktocpage=true,
	bookmarksnumbered=true,
	breaklinks=true,
	pdfstartview=FitH,
	hyperfigures=false,
	plainpages=false,
	naturalnames=true,
	colorlinks=true,
	pagebackref=false,
	pdfpagelabels]{hyperref}

\hypersetup{
	colorlinks,
	citecolor=blue,
	filecolor=blue,
	linkcolor=blue,
	urlcolor=blue
}

\usepackage[capitalize, noabbrev]{cleveref}
\let\subsectionSymbol\S
\crefname{subsection}{\subsectionSymbol\!\!}{subsections}

\calclayout

\makeatletter
\@namedef{subjclassname@2020}{%
	\textup{2020} Mathematics Subject Classification}
\makeatother

%% file: aux/usualcmds.tex
\newtheorem{theorem}[subsubsection]{Theorem}
\newtheorem*{theorem*}{Theorem}
\newtheorem{proposition}[subsubsection]{Proposition}
\newtheorem*{proposition*}{Proposition}

\newtheorem*{lemma*}{Lemma}

\newtheorem*{corollary*}{Corollary}

\theoremstyle{definition}
\newtheorem{definition}[subsubsection]{Definition}
\newtheorem*{definition*}{Definition}

\newtheorem*{remark*}{Remark}

\newtheorem*{example*}{Example}

\newtheorem*{construction*}{Construction}

\newtheorem*{convention*}{Convention}

\newtheorem*{terminology*}{Terminology}

\newtheorem*{notation*}{Notation}

\newtheorem*{question*}{Question}

\DeclareMathOperator{\face}{d}
\DeclareMathOperator{\dege}{s}
\DeclareMathOperator{\bd}{\partial}
\DeclareMathOperator{\sign}{sign}
\newcommand{\ot}{\otimes}
\DeclareMathOperator{\EZ}{EZ}

\newcommand{\N}{\mathbb{N}}
\newcommand{\Z}{\mathbb{Z}}

\newcommand{\R}{\mathbb{R}}
\renewcommand{\k}{\Bbbk}
\newcommand{\sym}{\mathbb{S}}

\newcommand{\gsimplex}{\mathbb{\Delta}}

\newcommand{\Fun}{\mathsf{Fun}}
\newcommand{\Set}{\mathsf{Set}}
\newcommand{\Top}{\mathsf{Top}}

\newcommand{\Ch}{\mathsf{Ch}}
\newcommand{\simplex}{\triangle}
\newcommand{\sSet}{\mathsf{sSet}}
\newcommand{\cube}{\square}

\DeclareMathOperator{\Sing}{Sing}

\DeclareMathOperator{\Sq}{Sq}

\DeclareMathOperator{\chains}{N}

\DeclareMathOperator{\gchains}{C}

\DeclarePairedDelimiter\bars{\lvert}{\rvert}

\DeclarePairedDelimiter\set{\{}{\}}

\newcommand{\id}{\mathsf{id}}
\renewcommand{\th}{\mathrm{th}}
\newcommand{\op}{\mathrm{op}}
\DeclareMathOperator*{\colim}{colim}

\newcommand{\xra}[1]{\xrightarrow{#1}}

\newcommand{\pdfEinfty}{\texorpdfstring{${E_\infty}$}{E-infty}}

\newcommand{\rD}{\mathrm{D}}

\newcommand{\rH}{\mathrm{H}}

\newcommand{\rS}{\mathrm{S}}

\newcommand{\cC}{\mathcal{C}}

\newcommand{\cE}{\mathcal{E}}

\newcommand{\cM}{\mathcal{M}}
\newcommand{\cN}{\mathcal{N}}

\newcommand{\cX}{\mathcal{X}}

\newcommand{\Sh}{\mathsf{Sh}}
\newcommand{\fI}{\mathfrak{I}}

\newcommand{\fX}{\mathfrak{X}}

%% file: aux/commands.tex
\DeclareMathOperator*{\hocolim}{hocolim}
\newcommand{\mSet}[1]{\mathsf{mSet}^{(#1)}}

\newcommand{\diag}{\rD}
\newcommand{\msimplex}[2]{\simplex^{#1}\times\dots\times\simplex^{#2}}
\newcommand{\gmsimplex}[2]{\gsimplex^{#1}\times\dots\times\gsimplex^{#2}}
\DeclareMathOperator{\gi}{\mathfrak{i}}

\newcommand{\UM}{\mathrm{U}\mathcal{M}}
\newcommand{\sur}{Sur}
\newcommand{\conf}[2]{Conf_{#1}(\R^{#2})}

\DeclareMathOperator{\schains}{S}
\DeclareMathOperator{\TR}{TR}
\DeclareMathOperator{\TC}{TC}
\DeclareMathOperator{\tc}{tc}
\newcommand{\BE}{E}
\newcommand{\CG}{\mathcal{K}}
\DeclareMathOperator{\rank}{rank}
\newcommand{\ex}{\mathrm{ex}}

\DeclareMathOperator{\ez}{\mathfrak{ez}}

\DeclareMathOperator{\In}{I}

\DeclareMathOperator{\copr}{\Delta}
\DeclareMathOperator{\aug}{\epsilon}
\DeclareMathOperator{\pr}{\ast}

%% file: sec/abstract.tex
% !TEX root = ../msimp.tex

\begin{abstract}
	This paper presents a generalization to multisimplicial sets of previously defined $E_\infty$-coalgebra structures on the chains of simplicial and cubical sets. We focus on the surjection chain complexes of McClure--Smith as a main example and construct a zig-zag of complexity preserving quasi-isomorphisms of $E_\infty$-coalgebras relating these to both the singular chains on configuration spaces and the Barratt--Eccles chain complexes.
\end{abstract}

%% file: sec/introduction.tex
% !TEX root = ../msimp.tex

\section{Introduction}\label{s:introduction}

The cochain complex of a simplicial set is equipped with the classical Alexander--Whitney product defining the ring structure in cohomology.
This cochain level structure has several explicit extensions to an $E_\infty$-algebra \cite{mcclure2003multivariable, berger2004combinatorial, medina2020prop1} encoding commutativity and associativity up to coherent homotopies.
The importance of $E_\infty$-algebras in homotopy theory is well known.
For example, Mandell showed that finite type nilpotent spaces are weakly equivalent if and only if their singular cochains are quasi-isomorphic as $E_\infty$-algebras \cite{mandell2006homotopy_type}.
Our first objective is to define a natural product together with an $E_\infty$-algebra extension on the cochains of multisimplicial sets \cite{gugenheim1957supercomplexes}.
These are generalizations of both simplicial and cubical sets which are useful for concrete computations since they can model homotopy types using fewer cells.
For example, the proof of the non-formality of the cochain algebra of planar configuration spaces \cite{salvatore2020planarnonformality} used a simplicial model and the Alexander--Whitney product on its cochains.
By using a multisimplicial model and the product defined here, these computations become simpler and faster, paving the way for extending this result to higher dimensions.

Multisimplicial sets are contravariant functors from products of the simplex category $\simplex$ to $\Set$.
Explicitly, for any positive integer $k$ the category $\mSet{k}$ of $k$-fold multisimplicial sets is the presheaf category $\Fun((\simplex^\op)^{\times k}, \Set)$.
There is a notion of geometric realization for multisimplicial sets, which results in a CW complex having, for each non-degenerate multisimplex, a cell modeled on a product of geometric simplices $\gmsimplex{n_1}{n_k}$.
We are interested in modeling homotopy types algebraically, for which we consider the composition of the geometric realization and the functor of cellular chains $\gchains$.
This composition defines $\chains \colon \mSet{k} \to \Ch$, the functor of (normalized) chains.
In \cref{ss:e-infty extension} we define a lift of $\chains$ to the category of $E_\infty$-coalgebras, and, consequently, a lift of the functor of cochains to the category of $E_\infty$-algebras.
We do so using the finitely presented $E_\infty$-prop introduced in \cite{medina2020prop1} and its monoidal properties.
Specifically, using the isomorphism
\[
\gchains(\gmsimplex{n_1}{n_k}) \cong
\gchains(\gsimplex^{n_1}) \ot\dotsb\ot \gchains(\gsimplex^{n_k}),
\]
we extend the image of the prop generators constructed in \cite{medina2020prop1} from the chains of standard simplices to those of standard multisimplices.
These generators are the Alexander--Whitney coproduct, the augmentation map, and an algebraic version of the join product.
The resulting $E_\infty$-coalgebra structure generalizes those defined in \cite{mcclure2003multivariable, berger2004combinatorial, medina2020prop1} for simplicial chains and in \cite{medina2022cube_einfty} for cubical chains.
As an application, we study the Steenrod construction for multisimplicial chains in \cref{ss:cup coproducts} emphasizing the explicit nature of our construction.

Let us now focus on the relationship between multisimplicial and simplicial theories.
The restriction to the image of the diagonal inclusion $\simplex^\op \to (\simplex^\op)^{\times k}$ of any $k$-fold multisimplicial set $X$ defines its associated diagonal simplicial set $X^\diag$.
There is a natural homeomorphism of realizations $\bars{X} \cong \bars{X^\diag}$
\cite{quillen2010higheralgktheory}.
Under this homeomorphism the cells of $\bars{X^\diag}$ arise from those of $\bars{X}$ through subdivision, a procedure described algebraically by the Eilenberg--Zilber quasi-isomorphism
\[
\EZ \colon \chains(X) \to \chains(X^\diag).
\]
The functor induced by the diagonal restriction has a right adjoint $\cN^{(k)}$, the multisimplicial nerve of a simplicial set, making the categories of $k$-fold multisimplicial and simplicial sets equivalent in Quillen's sense.
Furthermore, there is a natural inclusion
\[
\In \colon \chains(Y) \to \chains(\cN^{(k)} Y)
\]
which is also a quasi-isomorphism.
On one hand, the $\EZ$ map preserves the counital coalgebra structure, but it does not respect the higher $E_\infty$-structure.
On the other, the map $\In$ is an $E_\infty$-coalgebra quasi-isomorphism as proven in \cref{ss:inclusion}.
We use this fact to prove in \cref{ss:singular} that, for any topological space $\fX$, the linear map from its singular simplicial chains to its singular $k$-fold multisimplicial chains, given by precomposing a continuous map $(\gsimplex^n \to \fX)$ with the projection $(\gsimplex^n \times \gmsimplex{0}{0} \xra{\pi_1} \gsimplex^n)$, induces a natural quasi-isomorphism of $E_\infty$-coalgebras.

In the second part of the paper, we use these constructions to study a multisimplicial model of the canonical filtration
\[
\conf{r}{1} \subseteq \conf{r}{2} \subseteq \dotsb
\]
of the space $\conf{r}{\infty}$ of $r$ distinct ordered points in $\R^\infty = \colim (\R^1 \subseteq \R^2 \subseteq \dotsb)$.
Concretely, for any integer $r$, McClure and Smith \cite{mcclure2003multivariable} introduced a chain complex $\cX(r)$ of $\Z[\sym_r]$-modules with a filtration
\[
\cX_{1}(r) \subseteq \cX_{2}(r) \subseteq \dotsb
\]
and showed that $\cX(r)$ is connected to the singular chains of $\conf{r}{\infty}$ via a zig-zag of filtration preserving $\sym_r$-equivariant quasi-isomorphisms.
Presumably it was observed by both McClure--Smith and Berger--Fresse that $\cX(r)$ can be interpreted as the chains of an $r$-fold multisimplicial set $\sur(r)$, which we introduce in \cref{ss:surjection model} with a filtration
\[
\sur_{1}(r) \subseteq \sur_{2}(r) \subseteq \dotsb
\]
so that $\chains\sur_{d}(r) \cong \cX_{d}(r)$.
There is an operad structure on $\set{\cX_{d}(r)}_{r \geq 1}$ for each $d \geq 1$, but we do not focus on it since it is not induced from one at the multisimplicial level.
By the constructions in \cref{s:multisimplicial} the complex $\chains\sur(r)$ is equipped with an $E_\infty$-coalgebra structure, which we connect to the singular chains of $\conf{r}{\infty}$ via an explicit zig-zag of filtration preserving $\sym_r$-equivariant quasi-isomorphisms of $E_\infty$-coalgebras.

In a similar way, Berger and Fresse \cite{berger2004combinatorial} studied a chain complex $\cE(r)$ of $\Z[\sym_r]$-modules with a filtration
\[
\cE_{1}(r) \subseteq \cE_{2}(r) \subseteq \dotsb \,.
\]
This complex comes from the chains on a simplicial set introduced by Barratt and Eccles \cite{barrat1974operad} equipped with a filtration
\[
\BE_{1}(r) \subseteq \BE_{2}(r) \subseteq \dotsb
\]
due to Smith \cite{smith1989filtration}.
(As before we disregard the operadic structure.)
Since $\cE(r)$ is induced from a simplicial set, it is endowed with an $E_\infty$-coalgebra structure, and it is not hard to see that the zig-zag of filtration preserving $\sym_r$-equivariant quasi-isomorphisms used to compare it to the singular chains of $\conf{r}{\infty}$ respects this higher structure.
Consequently, $\cX(r)$ and $\cE(r)$ can be related by an explicit zig-zag of such maps.

It is desirable to have a direct map between the multisimplicial and simplicial models.
Berger--Fresse constructed two such filtrations preserving $\sym_r$-equivariant quasi-isomorphisms
\[
\TR \colon \chains\BE(r) \to \chains\sur(r)
\quad \text{and} \quad
\TC \colon \chains\sur(r) \to \chains \BE(r).
\]
The first one, introduced in \cite[1$\cdot$3]{berger2004combinatorial}, is unfortunately not a coalgebra map.
Therefore we will focus on the second one, which was introduced in \cite{berger2002prismatic}.
Our contribution, presented in \cref{ss:table completion}, is the construction of a factorization
\[
\TC \colon \chains\sur(r) \xra{\EZ} \chains\sur(r)^\diag \xra{\chains(\tc)} \chains \BE(r),
\]
where the second map is induced from a filtration preserving $\sym_r$-equivariant weak-equivalence of simplicial sets.
Therefore, we prove that $\TC$ is a coalgebra map since $\EZ$ is one.

%% file: sec/acknowledgment.tex
% !TEX root = ../msimp.tex

\section*{Acknowledgments}

Some of the results of this work appeared in the B.Sc. thesis of A.P. written under the guidance of P.S.: `Generalization of Alexander--Whitney map to Multisimplicial Sets', University of Rome `Tor Vergata' (2019).

A.M-M. would like to thank Peter May and Dennis Sullivan for insightful comments about this project.

A.P. and P.S. were partially supported by the MIUR Excellence Department Project MatMod@TOV awarded to the Department of Mathematics, University of Rome `Tor Vergata'.

A.M-M. Would like the acknowledge the support and excellent working conditions of the Max Planck Institute for Mathematics in Bonn.

%% file: sec/multisimplicial.tex
% !TEX root = ../msimp.tex

\section{Multisimplicial algebraic topology}\label{s:multisimplicial}

\subsection{Multisimplicial sets}\label{ss:multisimplicial sets}

Let us consider an arbitrary positive integer~$k$.
The \textit{$k$-fold multisimplex category} $\simplex^{\times k}$ is the $k$-fold Cartesian product of the simplex category $\simplex$.
The category
\[
\mSet{k} = \Fun\big((\simplex^{\times k})^\op, \Set\big)
\]
is referred to as the category of \textit{$k$-fold multisimplicial sets}.
We remark that $\mSet{1}$ and $\mSet{2}$ are naturally equivalent to the categories of simplicial and bisimplicial sets respectively.
A representable $k$-fold multisimplicial sets is denoted by $\simplex^{n_1,\dots,n_k}$.

Explicitly, a $k$-fold multisimplicial set $X$ consists of a collection of sets
\[
X_{m_1,\dots,m_k} = X \big( [m_1] \times\dots\times [m_k] \big)
\]
indexed by $k$-tuples of non-negative integers $(m_1,\dots,m_k)$ together with \textit{face maps}
\[
\face_i^j \colon X_{m_1, \dots, m_j,\dots,m_k} \to X_{m_1, \dots, m_j-1, \dots, m_k}
\]
and \textit{degeneracy map}
\[
\dege^j_i \colon X_{m_1, \dots, m_j,\dots,m_k} \to X_{m_1, \dots, m_j+1, \dots, m_k}
\]
for $1 \leq j \leq k$ and $0 \leq i \leq m_j$ such that, referring to $j$ as the \textit{direction} of these maps, two of them satisfy the simplicial identities when they have the same direction and commute when they do not.
An element of $X_{m_1,\dots,m_k}$ is called an $(m_1,\dots,m_k)$-\textit{multisimplex} and it is said to be \textit{degenerate} if it is in the image of a degeneracy map.

\subsection{Geometric realization}\label{ss:geometric realization}

We will use the following model of the topological simplex:
\[
\gsimplex^n = \big\{
(t_1, \dots, t_n) \in [0,1]^n \mid t_1 \geq \dots \geq t_n
\big\}
\]
with
\[
\delta_i(t_1, \dots, t_n) =
\begin{cases}
	(1, t_1, \dots, t_n) & i = 0, \\
	(t_1, \dots, t_i, t_i, \dots, t_n) & 0 < i < n, \\
	(t_1, \dots, t_n, 0) & i = n,
\end{cases}
\]
and
\[
\sigma_i(t_1, \dots, t_n) = (t_1, \dots, \widehat t_i, \dots, t_n).
\]

The \textit{geometric realization} functor
\[
\bars{-} \colon \mSet{k} \to \Top
\]
is the Yoneda extension of the functor defined on representable objects by
\[
\bars{\simplex^{n_1,\dots,n_k}} = \gmsimplex{n_1}{n_k}.
\]

Explicitly, for a $k$-fold multisimplicial set $X$ we have
\[
\bars{X} \cong
\coprod \gmsimplex{n_1}{n_k} \times X_{n_1,\dots,n_k} \, /_\sim
\]
where
\[
\begin{split}
	(\vec{t}_1, \dots, \vec{t}_j, \dots, \vec{t}_k, \face^j_i(x)) &\sim (\vec{t}_1, \dots, \delta_i(\vec{t}_j), \dots, \vec{t}_k, x), \\
	(\vec{t}_1, \dots, \vec{t}_j, \dots, \vec{t}_k, \dege^j_i(x)) &\sim (\vec{t}_1, \dots, \sigma_i(\vec{t}_j), \dots, \vec{t}_k, x),
\end{split}
\]
which equips $\bars{X}$ with a canonical cellular structure.

The geometric realization functor has a right adjoint
\[
\Sing^{(k)} \colon \Top \to \mSet{k}
\]
defined on a topological space $\fX$, as usual, by the expression
\[
\Sing^{(k)}(\fX)_{n_1,\dots,n_k} =
\Top(\gmsimplex{n_1}{n_k},\fX).
\]

\subsection{Algebraic realization}\label{ss:algebraic realization}

The functor of \textit{chains}
\[
\chains \colon \mSet{k} \to \Ch,
\]
is the Yoneda extension of the functor defined on representable objects by
\[
\chains \big( \simplex^{n_1, \dots, n_k} \big) =
\chains(\simplex^{n_1}) \ot \dotsb \ot \chains(\simplex^{n_k}).
\]
It is naturally isomorphic to the composition of the geometric realization functor and the functor of cellular chains with respect to the canonical cellular structure.

Explicitly, for a $k$-fold multisimplicial set $X$ the $\k$-module $\chains(X)_n$ is freely generated by the non-degenerate $(n_1, \dots, n_k)$-multisimplices with $n_1+\dots+n_k = n$.
The differential $\bd \colon \chains(X)_n \to \chains(X)_{n-1}$ is given on one such basis element by
\[
\bd(x) = \sum_{j=1}^k \sum_{\ell_j=1}^{n_j}
(-1)^{n_{1}+\dots+n_{j-1}+\ell_j} \, \face^j_{\ell_j}(x).
\]

For any topological space $\fX$ the chain complex $\chains\Sing^{(k)}(\fX)$ is denoted $\rS^{(k)}(\fX)$ and referred to as the $k$-\textit{fold singular chains} of $\fX$.

\subsection{Coalgebra structure}\label{ss:coalgebra}

A \textit{counital coalgebra} structure on a chain complex $C$ is a pair of chain maps $\copr \colon C \to C \ot C$ and $\aug \colon C \to \k$ satisfying
\[
(\id \ot \aug) \circ \copr =
\id =
(\aug \ot\, \id) \circ \copr.
\]
The tensor product of two counital coalgebras $C$ and $C'$ is itself a counital coalgebra with structure maps given by
\begin{gather*}
	C \ot C^\prime \xra{\copr \ot \copr^\prime}
	(C \ot C) \ot (C^\prime \ot C^\prime) \xra{\tau}
	(C \ot C^\prime) \ot (C \ot C^\prime), \\
	C \ot C^\prime \xra{\aug \ot \aug^\prime}
	\k \ot \k \xra{\cong} \k,
\end{gather*}
where $\tau$ transposes the second and third factors.

For each $n \in \N$, the complex $\chains(\simplex^n)$ is naturally equipped with a counital coalgebra structure defined by:
\[
\begin{split}
	\copr \big( [v_0, \dots, v_m] \big) &=
	\sum_{i=0}^m \, [v_0, \dots, v_i] \ot [v_i, \dots, v_m], \\
	\aug \big( [v_0, \dots, v_q] \big) &=
	\begin{cases} 1 & \text{ if } q = 0, \\ 0 & \text{ if } q>0. \end{cases}
\end{split}
\]
We will refer to it as the \textit{Alexander--Whitney structure}.

Using the tensor product structure, we deduce a natural counital coalgebra structure on the chains of representable multisimplicial sets
\[
\chains \big( \simplex^{n_1, \dots, n_k} \big) =
\chains(\simplex^{n_1}) \ot \dotsb \ot \chains(\simplex^{n_k})
\]
and, via a Yoneda extension, one on the chains of general multisimplicial sets.

Explicitly, for a $k$-fold multisimplicial set $X$ and $(m_1,\dots,m_k)$-multisimplex $x$ let
\[
\fI_{m_1,\dots,m_k} = \set[\big]{(i_{1},\dots,i_{k}) \mid 0 \le i_j \le m_j,\ \forall j = 1,\dots,k},
\]
then
\[
\copr(x) =
\sum_{I\in \mathfrak{I}_{k,x}} \;
(-1)^{\sum_{1 \leq l<h \leq k} i_h (m_l-i_l) } \
x \rfloor_{(i_{1},\dots,i_{k})} \ot
\!\,_{(m_{1}-i_{1}, \dots, m_{k}-i_{k})} \lfloor x
\]
where the \textit{front $(i_1,\dots,i_k)$-face} of $x$ is the multisimplex
\[
x \rfloor_{(i_{1}, \dots, i_{k})} =
X(F_{i_1}, \dots, F_{i_k})(x) \in X_{i_1,\dots,i_k}
\]
with
$F_{i_j} \colon [i_j] \to [n_j]$ defined by $F_{i_j}(h)=h$, and the \textit{back $(i_1,\dots,i_k)$-face} of $x$ is the multisimplex
\[
\,_{(i_{1}, \dots, i_{k})} \lfloor x =
X(B_{i_1}, \dots, B_{i_k})(x) \in X_{i_1,\dots,i_k}
\]
with $B_{i_j} \colon [i_j] \to [n_j]$ defined by $B_j(h) = h+m_j-i_j$.

\subsection{\pdfEinfty-extension}\label{ss:e-infty extension}

An \textit{$\cM$-bialgebra} is a counital coalgebra $(C, \copr, \aug)$ together with a degree $1$ linear map $\pr \colon C \ot C \to C$ satisfying
\[
\bd (c_1 \ast c_2) - \bd c_1 \ast c_2 + (-1)^{\bars{c_1}} c_1 \ast \bd c_2 =
\aug(c_1) c_2 - \aug(c_2) c_1,
\]
\[
\aug (c_1 \pr c_2) = 0,
\]
for all $c_1, c_2 \in C$.
As proven in \cite{medina2020prop1}, the collection of all maps $\set{C \to C^{\ot r}}_{r \in \N}$ generated by $\copr$, $\aug$ and $\pr$ make $C$ into an $E_\infty$-coalgebra, that is to say, a coalgebra over certain operad $\UM$ that is a cofibrant resolution of the terminal operad.

As proven in \cite{medina2021cobar}, the counital coalgebra structure on the tensor product of two $\cM$-bialgebras $C$ and $C'$ can be naturally extended to an $\cM$-bialgebra structure using
\[
(C \ot C^\prime) \ot (C \ot C^\prime) \xra{\tau}
C \ot C \ot C^\prime \ot C^\prime
\xra{\ \aug \ot \, \id \, \ot \, \pr + \pr \ot \, \id \, \ot \,\aug\ }
C \ot C^\prime.
\]

For any integer $n$, the \textit{join product} $\ast \colon \chains(\simplex^n)^{\ot 2} \to \chains(\simplex^n)$ is the natural degree~$1$ linear map defined by
\begin{equation*}
	\left[v_0, \dots, v_p \right] \pr \left[v_{p+1}, \dots, v_q\right] =
	\begin{cases} (-1)^{p} \sign(\pi) \left[v_{\pi(0)}, \dots, v_{\pi(q)}\right] & \text{ if } v_i \neq v_j \text{ for } i \neq j, \\
		\hfil 0 & \text{ if not}, \end{cases}
\end{equation*}
where $\pi$ is the permutation that orders the vertices.
It is proven in \cite{medina2020prop1} that on the chains of representable simplicial sets the Alexander--Whitney structure together with the join product make $\chains(\simplex^n)$ into a natural $\cM$-bialgebra and, consequently, a natural $E_\infty$-coalgebra.
We mention that this structure is induced by one preset at the level of geometric realizations \cite{medina2021prop2}.

Using the tensor product structure, we deduce a natural $\cM$-bialgebra structure on the chains of representable multisimplicial sets
\[
\chains\big(\simplex^{n_1,\dots,n_k}\big) =
\chains(\simplex^{n_1}) \ot\dotsb\ot \chains(\simplex^{n_k}),
\]
and consequently a natural $E_\infty$-coalgebra structure, which extends along the Yoneda inclusion to the chains on any multisimplicial set $X$.

Explicitly, for two basis elements of $\chains\big(\simplex^{n_1,\dots,n_k}\big)$ we have
\[
(x_1 \ot \dotsb \ot x_n) \ast (y_1 \ot \dotsb \ot y_n) =
\sum_{i=1}^n x_{<i}\, \epsilon(y_{<i}) \ot x_i \ast y_i \ot \epsilon(x_{>i}) \, y_{>i},
\]
where, with the convention $x_{<1} = x_{>n} = 1 \in \k$,
\begin{align*}
	z_{<i} & = z_1 \ot \dotsb \ot z_{i-1}, \\
	z_{>i} & = z_{i+1} \ot \dotsb \ot z_n.
\end{align*}

We remark that since the category of $\cM$-bialgebras is not cocomplete, we do not necessarily have an $\cM$-bialgebra structure on $\chains(X)$ for a general multisimplicial set $X$.
An example for which such structure does not exist is given by one such $X$ whose geometric realization consists of just two points.

\subsection{Cubical theory}\label{ss:cubical}

Since the complex of chains of the $k$-fold multisimplicial set $\msimplex{1}{1}$ is isomorphic to the chains on the standard cubical set $\cube^k$, it is natural to compare the $E_\infty$-coalgebra structure defined here with that presented in \cite{medina2022cube_einfty} for cubical sets.
As counital coalgebras $\chains(\msimplex{1}{1})$ and $\chains(\cube^k)$ are isomorphic, and, denoting the product of the $\cM$-bialgebra defined there by $\widetilde\ast$, we have
\[
x \ \widetilde\ast \ y = (-1)^{\bars{x}} x \ast y
\]
under this chain isomorphism.
The sign convention used here is more natural, used for example to endow Adams' cobar construction with the structure of a monoidal $E_\infty$-coalgebra \cite{medina2021cobar}.

\subsection{Steenrod construction}\label{ss:cup coproducts}

In \cite{steenrod1947products}, Steenrod introduced natural operations on the mod~2 cohomology of spaces, the celebrated \textit{Steenrod squares}
\[
\begin{tikzcd} [column sep=small, row sep=0]
	\Sq^k \colon &[-20] \rH^{-n} \arrow[r] & \rH^{-n-k} \\ &
	{[\alpha]} \arrow[r, mapsto] & \big[ (\alpha \ot \alpha) \Delta_{n-k} \big],
\end{tikzcd}
\]
via an explicit construction of natural linear maps $\Delta_i \colon \chains(X) \to \chains(X) \ot \chains(X)$ for any simplicial set $X$, satisfying up to signs the following homological relations
\[
\bd \circ \, \Delta_i + \Delta_i \circ \bd = (1 + T) \Delta_{i-1},
\]
with the convention $\Delta_{-1} = 0$.
These so-called \textit{cup-$i$ coproducts} appear to be fundamental, as they are axiomatically characterized \cite{medina2022axiomatic} and induce the nerve of strict infinity categories \cite{medina2020globular}.
A description of cup-$i$ coproducts for multisimplicial sets can be deduced from our $E_\infty$-coalgebra structure.
It is given recursively by
\[
\begin{split}
	& \Delta_0 = \Delta, \\
	& \Delta_i =
	(\ast \ot \id) \circ (23) \circ (\Delta_{i-1} \ot \id) \circ \Delta.
\end{split}
\]

Steenrod also introduced operations on the mod~$p$ cohomology of spaces when $p$ is an odd prime \cite{steenrod1952reduced, steenrod1953cyclic}.
To define these effectively, generalization of the cup-$i$ coproducts were introduced in \cite{medina2021may_st}.
After the present work, these so-called \textit{cup-$(p,i)$ coproducts} are defined on multisimplicial chains, and their formulas are explicit enough to be implemented in the computer algebra system \href{https://comch.readthedocs.io/en/latest/}{\texttt{ComCH}} \cite{medina2021comch}, where constructions of Cartan and Adem coboundaries \cite{medina2020cartan,medina2021adem} for multisimplicial sets are also be found.

%% file: sec/comparison.tex
% !TEX root = ../msimp.tex

\section{Comparison with the simplicial theory}\label{s:comparison}

We will use $\sSet$ to denote the category of $1$-fold multisimplicial sets $\mSet{1}$ referring to its objects and morphisms as simplicial sets and simplicial morphisms as usual.

\subsection{Diagonal simplicial set}\label{ss:diagonal}

For any $k \in \N$, the diagonal
\[
\simplex^\op \xra{\diag}
(\simplex^\op)^{\times k} \xra{\cong}
(\simplex^{\times k})^{\op}
\]
induces a functor
\[
(-)^{\diag} \colon \mSet{k} \to \sSet
\]
explicitly defined on a $k$-fold multisimplicial set $X$ by
\begin{gather*}
	X^{\diag}_m = X_{m, \dots, m},
	\qquad
	\face_i = \face_i^1 \circ \dots \circ \face_i^k,
	\qquad
	\dege_i = \dege_i^1 \circ \dots \circ \dege_i^k.
\end{gather*}
It is straightforward to verify that
\[
\big(\simplex^{n_1,\dots,n_k}\big)^\diag \cong
\msimplex{n_1}{n_k}
\]
as simplicial sets.

The functor $(-)^\diag \colon \mSet{k} \to \sSet$ admits a right adjoint $\cN^{(k)} \colon \sSet \to \mSet{k}$, defined, as usual, by the expression
\[
\cN^{(k)}(Y)_{m_1,\dots,m_k} =
\sSet\big( \simplex^{m_1} \times \dots \times \simplex^{m_k}, \, Y \big).
\]
These functors define a Quillen equivalence.
A proof of this fact can be given using \cite[Proposition~1.6.8]{maltsiniotis2005grothendieck} or adapting that in \cite[Proposition~1.2]{moerdijk1989bisimplicialsets}.

\subsection{Eilenberg--Zilber map}\label{ss:eilenber-zilber}

Recall that an \textit{$(n_1,\dots,n_k)$-shuffle} $\sigma$ is a permutation in $\sym_n$ satisfying
\begin{gather*}
	\sigma(1) < \dots < \sigma(n_1), \\
	\sigma(n_1+1) < \dots < \sigma(n_1+n_2), \\
	\vdots \\
	\sigma(n-n_k-1) < \dots < \sigma(n),
\end{gather*}
where $n = n_1+\dots+n_k$.
We denote the set of such permutations by $\Sh(n_1,\dots,n_k)$.
For any $\sigma \in \Sh(n_1,\dots,n_k)$ the inclusion
\[
\gi_\sigma \colon \gsimplex^n \to \gmsimplex{n_1}{n_k}
\]
is defined by the assignment
\[
(x_1,\dots,x_n) \mapsto (x_{\sigma^{-1}(1)}, \dots, x_{\sigma^{-1}(n)}).
\]
If $e$ is the identity permutation, we denote $\gi_{e}$ simply as $\mathfrak{i}$.
The set $\set{\gi_\sigma \mid \sigma \in \Sh(n_1,\dots,n_k)}$ defines a triangulation of $\gmsimplex{n_1}{n_k}$ making it isomorphic, in the category of cellular spaces, to the geometric realization of the simplicial set $\msimplex{n_1}{n_k}$.
Using this identification, the identity map induces a cellular map
\[
\ez \colon \gmsimplex{n_1}{n_k} \to \bars[\big]{\msimplex{n_1}{n_k}},
\]
whose induced chain map
\[
\EZ \colon \chains(\simplex^{n_1,\dots,n_k}) \to \chains\big(\msimplex{n_1}{n_k}\big),
\]
agrees, under the natural identifications, with the traditional Eilenberg--Zilber map.

For a multisimplicial set $X$, the induced chain map $\EZ \colon \chains(X) \to \chains(X^\diag)$ is explicitly given on an $(n_1,\dots,n_k)$-multisimplex $x$ by
\[
\EZ(x) \ = \sum_{\qquad \crampedclap{\sigma \in \Sh(n_1,\dots,n_k)}} \sign(\sigma) \, X(\sigma_1 ,\dots, \sigma_k)(x)
\]
where, for $\ell \in \set{1,\dots,k}$, the morphisms $\sigma_\ell \colon [n] \to [n_\ell]$ are defined by the following property: for
each $j \in \set{1,\dots,n}$ there is exactly one $\ell \in \set{1,\dots,k}$ such that $\sigma_\ell(j+1) = \sigma_\ell(j)+1$ and $\sigma_i(j+1) = \sigma_i(j)$ for all $i \neq \ell$.

Since the traditional Eilenberg--Zilber map preserves counital coalgebra structures we have the following.

\begin{theorem}
	For every multisimplicial set $X$ the map $\EZ \colon \chains(X) \to \chains(X^\diag)$ is a quasi-isomorphism of counital coalgebras.
\end{theorem}

We remark that the Eilenberg--Zilber map is not a morphism of $E_\infty$-coalgebras.
For example, as shown in \cite[\S5.4]{medina2022cube_einfty}, we have
\[
\Delta_1 \circ \EZ\big([0,1] \ot [0,1]\big) \neq
\EZ ^{\otimes 2}\, \circ \, \Delta_1\big([0,1] \ot [0,1]\big).
\]

\subsection{Canonical inclusion}\label{ss:inclusion}

Let $Y$ be a simplicial set and $n$ an integer.
Consider the function $Y_n \to \cN^{(k)} Y_{n,0,\dots,0}$ sending a simplex with characteristic map $\zeta \colon \simplex^n \to Y$ to the composition
\[
\simplex^n \times \msimplex{0}{0} \xra{\pi_1} \simplex^n \xra{\zeta_y} Y.
\]
These functions induce a chain map
\[
\In \colon \chains(Y) \to \chains(\cN^{(k)} Y)
\]
and we have the following.

\begin{theorem}
	The canonical inclusion $\In \colon \chains(Y) \to \chains(\cN^{(k)} Y)$ is a quasi-isomorphism of $E_\infty$-coalgebras for any simplicial set $Y$.
\end{theorem}

\begin{proof}
	The structure-preserving properties of this map are immediate.
	It remains to be shown that it induces a homology isomorphism.
	Consider the composition of quasi-isomorphisms
	\[
	\chains(\cN^{(k)} Y) \xra{\EZ}
	\chains\big((\cN^{(k)} Y)^\diag\big) \to
	\chains(Y)
	\]
	where the second map is induced by the counit of the adjunction.
	We will now verify that it is left inverse to $\In$.
	Consider a simplex $y$ with characteristic map $\zeta \colon \simplex^n \to Y$.
	The multisimplex $\In(y)$ is given by the simplicial map $\simplex^n \times \msimplex{0}{0} \xra{\pi_1} \simplex^n \xra{\zeta} Y$.
	Since the only $(n,0,\dots,0)$-shuffle is the identity, the simplex $\EZ \circ \In (y)$ is the simplicial map
	\[
	\zeta \circ \pi_1 \colon \simplex^n \times \msimplex{n}{n} \to Y.
	\]
	Finally, the image of this simplex under the counit is the evaluation of $[n] \times\dots\times [n]$ on $\zeta \circ \pi_1$ which gives $\zeta[n] = y$ as claimed.
\end{proof}

\subsection{Singular chains}\label{ss:singular}

\begin{theorem}
	Let $\fX$ be a topological space.
	The chain map
	\[
	\rS(\fX) \to \rS^{(k)}(\fX),
	\]
	defined by precomposing a continuous map $(\gsimplex^n \to \fX)$ with the projection
	\[
	\gsimplex^n \times \gmsimplex{0}{0} \xra{\pi_1} \gsimplex^n,
	\]
	is a quasi-isomorphism of $E_\infty$-coalgebras.
\end{theorem}

\begin{proof}
	This map factors as the composition of two quasi-isomorphisms of $E_\infty$-coalgebras.
	The first is $\In \colon \rS(\fX) \to \chains(\cN^{(k)} \Sing(\fX))$, which was studied in \cref{ss:inclusion}.
	The second is induced by a multisimplicial isomorphism
	\[
	\cN^{(k)} \Sing(\fX) \to \Sing^{(k)}(\fX)
	\]
	defined as follows.
	Using the adjunction of \cref{ss:geometric realization}, any simplicial map $\msimplex{n_1}{n_k} \to \Sing(\fX)$ corresponds canonically to a continuous map $\bars{\msimplex{n_1}{n_k}} \to \fX$, which precomposing with $\ez$ gives a continuous map $\gmsimplex{n_1}{n_k} \to \fX$.
	It is not hard to see that every such map arises this way since $\ez$ is a homeomorphism.
\end{proof}

%% file: sec/configuration.tex
% !TEX root = ../msimp.tex

\section{Models of configuration spaces}

We are interested in modeling algebraically the $\sym_r$-equivariant homotopy type of the space of configurations of $r$ labeled and distinct points in Euclidean $d$-dimensional space.
Multisimplicial sets can be used to provide an explicit chain complex model with a small number of generators, which, using the $E_\infty$-structure defined in this paper, retains all homotopical information by Mandell's theorem \cite{mandell2006homotopy_type}.

In the first subsection, we recall a method due to Berger detecting spaces homotopy equivalent to euclidean configuration spaces by means of a filtration indexed by a {\em complete graph poset}.
In the second subsection, we construct the multisimplicial model and show that is equipped with such a filtration.
In the third subsection, we recall the construction of the simplicial Barratt--Eccles model and show that is equipped with a similar filtration.
In the fourth subsection, we relate the multisimplicial and simplicial chain models by an explicit map. In the last subsection, we give some examples of the sizes of the two models, showing that the multisimplicial is smaller.

\input{sec/recognition}
\input{sec/surjections}
\input{sec/barratt-eccles}

\input{sec/tc}
\input{sec/count}

%% file: sec/recognition.tex
% !TEX root = ../msimp.tex

\subsection{Recognition of configuration spaces}\label{ss:recognition}

Let $\conf{r}{d}$ denote the configuration space of $r$-tuples of pairwise disjoint vectors in $\R^{d}$.
This space is equipped with a free action of the symmetric group $\sym_r$ of permutations of $\{1,\dots,r\}$ swapping elements of a $r$-tuple.

\begin{definition}
	A \textit{complete graph} on $r$ vertices is a pair $(\mu,\sigma)$ with $\mu$ a collection of non-negative integers $\mu_{ij}$ for all $1 \leq i < j \leq r$, and $\sigma$ is an ordering of $\{1,\dots,r\}$.
	We write $\sigma_{ij}$ for the restriction of the ordering $\sigma$ to the set $\{i,j\}$.
	Graphically $(\mu,\sigma)$ is a simple directed graph in the edge corresponding to $i<j$ directed according to $\sigma_{ij}$ and labeled by $\mu_{ij}$.
	Please consult \cref{f:complete graph} for an example.
	Let us denote the set of complete graphs with $r$ vertices by $\CG(r)$ equipped with the poset structure
	\begin{equation*}
		(\mu,\sigma)\le (\nu,\tau) \ \ \iff \ \
		\forall i,j \ (\mu_{ij}<\nu_{ij}) \ \text{ or } \
		(\mu_{ij},\sigma_{ij})= (\nu_{ij},\tau_{ij})
	\end{equation*}
	for each pair $i<j$.
	It is equipped with an exhaustive filtration by subposets
	\[
	\CG_1(r) \subset \CG_2(r) \subset \dotsb
	\]
	where $\CG_d(r)$ consists of those graphs with $\max(\mu_{ij}) < d$.
\end{definition}

\begin{figure}
	\centering
	\begin{equation*}
		\begin{tikzcd}
			\ & 2 & \ \\
			1 \arrow[ur,"2"] \arrow[rr] \arrow[dr,"3"']& \arrow[r,"1"] \arrow[u,"2"] & 3 \arrow[ul,"1"'] \\
			\ & 4 \arrow[ur,"4"'] \arrow[uu] & \
		\end{tikzcd}
	\end{equation*}
	\caption{A complete graph on 4 vertices with ordering $\sigma=(1432)$ and $\mu=(\mu_{12},\mu_{13},\mu_{14},\mu_{23},\mu_{24},\mu_{34})=(2,1,3,1,2,4)$.}
	\label{f:complete graph}
\end{figure}

\begin{definition}\label{cellulardecomposition}
	For a given poset $A$, a cellular $A$-decomposition of a topological space $\fX$ is a family of subspaces $\set{\fX_a}_{a \in A}$ such that:
	\begin{itemize}
		\item [i.] $a \leq b$ implies $\fX_a \subseteq \fX_b$;
		\item [ii.] $\colim_{a \in A} \fX_a = \fX$;
		\item [iii.] $\fX_a$ is contractible for each $a$;
		\item [iv.] $\bigcup_{a<b} \fX_a \subset \fX_b$ is a closed cofibration.
	\end{itemize}
\end{definition}

The relevance of this notion is the well-known fact that if a topological space $\fX$ admits a \textit{cellular $A$-decomposition}, then the natural maps
\begin{equation}\label{e:cellular poset decomposition}
	\fX = \colim_A \fX_a \leftarrow \hocolim_A \fX_a \to \bars{A}
\end{equation}
are cellular homotopy equivalences.
Please consult \cite[\S1.7]{berger1997confspacemodel} for a proof.

Let $\cC_d(r)$ be the space of $r$ little $d$-dimensional cubes, which is equipped with an equivariant homotopy equivalence to $\conf{r}{d}$ picking the center of cubes.
Brun and others in \cite{brunfiedorowiczvogt2007} show that $\cC_d(r)$ has a cellular $\CG^\ex_d(r)$-decomposition $\{\cC_a\}$, where $\CG^\ex_d(r)$ is a poset containing the poset $\CG_d(r)$ and the inclusion of posets induces an equivariant homotopy equivalence on realizations.
Combining these results we have

\begin{proposition}\label{p:zig-zag conf}
	If a space $\fX$ has a cellular $\CG_d(r)$-decomposition, then $\fX$ is equal to $\colim_{\CG_d(r)} \fX_a$ and
	\begin{equation*}
		\begin{tikzcd}[column sep=small, row sep=5]
			\displaystyle \colim_{\CG_d(r)} \fX_a & \arrow[l] \displaystyle \hocolim_{\CG_d(r)} \fX_a \arrow[r] & \bars{\CG_d(r)} \arrow[d] & & \\ & &
			\bars{\CG^\ex_d(r)} & \arrow[l] \displaystyle \hocolim_{\CG^\ex_d(r)}\,\cC_\alpha \arrow[r] & \cC_d(r) \arrow[r] & \conf{r}{d}
		\end{tikzcd}
	\end{equation*}
	is a zig-zag of equivariant homotopy equivalences.
\end{proposition}

\begin{definition}
	Let $X$ be a multisimplicial (or simplicial) set.
	A \textit{$\CG(r)$-filtration} of $X$ is a family of (multi)simplicial subsets $\{X_a\}$ indexed by
	$a \in \CG(r)$ so that
	\begin{enumerate}
		\item $a \leq b$ implies $X_a \subseteq X_b$;
		\item $|X_a|$ is a cellular $\CG(r)$-decomposition of the realization $|X|$
	\end{enumerate}
	In particular this implies that $X=\colim_{a \in \CG(r)}X_a$.
	Let $X_d = \colim_{a \in CG_d(r)} X_a$.
	There is a nested sequence
	\[
	X_1 \subset X_2 \subset \dotsb
	\]
	For a given (multi)simplex $x \in X$ we will refer to $\min\set{d \mid x \in X_d}$ as the \textit{complexity} of $x$.
\end{definition}

%% file: sec/surjections.tex
% !TEX root = ../msimp.tex

\subsection{Multisimplicial model}\label{ss:surjection model}

We define for each positive integer $r$ a multisimplicial set $\sur(r)$ equipped with a $\CG(r)$-filtration.
The functor of chains applied to the nested sequence
\[
\sur_{1}(r) \subset \sur_{2}(r) \subset \dotsb
\]
 will recover the algebraic models
\[
\chi_1(r) \subset \chi_2(r) \subset \dotsb
\]
of configuration spaces developed by McClure--Smith \cite{mcluresmith2004geomodel}.

Spaces $Y_0^r$ homeomorphic to $\bars{\sur(r)}$ were studied in the work of McClure--Smith \cite{mcclure2003multivariable}.
The homeomorphism between $Y_0^r$ and $\bars{\sur(r)}$ is described explicitly in the appendix of \cite{salvatore2009deligne}.

Let $\sur(r)$ be the $k$-fold multisimplicial set that has as $(m_1,\dots,m_r)$-multisimplices the surjective maps
\[
f \colon \set{1,\dots,m+r} \to \set{1,\dots,r},
\]
where $m = m_1+\dots+m_r$, satisfying that the cardinality of $f^{-1}(\ell)$ is $m_\ell$ for each $\ell \in \set{1,\dots,r}$.
We represent this multisimplex by the sequence $f(1) \dotsm f(m+r)$.
The face and degeneracy maps $\face^\ell_j$ and $\dege^\ell_j$ act on it by respectively removing and doubling the $(j+1)^\th$ occurrence of $\ell$ in the sequence.

Next we define a $\CG(r)$-filtration on $\sur(r)$.
For $i<j$, let $f_{ij}$ be the subsequence of $f(1) \dotsm f(m+r)$ obtained by omitting all occurrences of elements different from $i$ and $j$.
For $(\mu,\sigma) \in \CG(r)$
we say that 
$f \in \sur(r)_{(\mu,\sigma)}$
 if for each $i<j$, either $i$ and $j$ alternate strictly less than $\mu_{ij}$ times in the sequence $f_{ij}$, or they do so exactly $\mu_{ij}$ times and the ordering formed by the first occurrences of $i$ and $j$ in $f_{ij}$ agrees with $\sigma_{ij}$.
%reference of proof

The surjection $f$ has complexity $d$ or less if the alternation number of each $f_{ij}$ is less than $d+1$, i.e., if the non-degenerate dimension of $f_{ij}$ in $\sur(2)$ is $d$ or less for each $i<j$.
We notice that the action of $\sym_r$ on $\sur(r)$ preserves the nested sequence
$$\sur_1(r) \subset \sur_2(r) \subset \dots$$
For the proof that $|\sur(r)|$
has indeed an induced cellular $\CG(r)$-decomposition we refer to 
Lemma 14.8 in \cite{mcluresmith2004geomodel}.
Applying the functor of singular chains to the zig-zag of \cref{p:zig-zag conf} produces a zig-zag of equivariant quasi-isomorphisms of $\UM$-coalgebras connecting $\schains\bars{\sur_d(r)}$ and $\schains\conf{r}{d}$.
We can extend it using the following zig-zag of maps of the same kind
\[
\schains\bars{\sur_{d}(r)} \cong
\schains\bars{\sur_{d}(r)^{D}} \to
\chains(\sur_{d}(r)^D) \to
\chains(\cN^{(r)}(\sur_{d}(r)^D)) \leftarrow
\chains\sur_{d}(r).
\]
The first map is induced by the homeomorphism $\bars{\sur_{d}(r)} \cong \bars{\sur_{d}(r)^{D}}$, the second by the unit of the Quillen equivalence between simplicial sets and topological spaces, the third is the comparison map of \cref{ss:inclusion}, and last one is induced by the unit of the Quillen equivalence between multisimplicial sets and simplicial sets.

As announced in the introduction, this construction relates the chains on the multisimplicial model of configuration space and its singular chains via an explicit zig-zag of equivariant quasi-isomorphisms of $E_\infty$-coalgebras.

%% file: sec/barratt-eccles.tex
% !TEX root = ../msimp.tex

\subsection{Simplicial model}\label{ss:simplicial model}

We recall the Barratt--Eccles simplicial set $\BE(r)$ defined for each $r\in\N$ that is equipped with a $\CG(r)$-filtration.
Applying the functor of chains to the nested sequence
\[
\BE_{1}(r) \subset \BE_{2}(r) \subset \dotsb.
\]
will provide the algebraic models
\[
\cE_{1}(r) \subset \cE_{2}(r) \subset \dotsb
\]
of configuration spaces studied by Berger and Fresse in \cite{berger2004combinatorial}.

The $n$-simplices of $\BE(r)$ are tuples of $n+1$ elements of the symmetric group $\sym_r$.
Its face and degeneracy maps are defined by removing and doubling elements respectively.
There is an operad structure on these simplicial sets, but we do not consider it here.

Next we recall a $\CG(r)$-filtration on $\BE(r)$.
For $i<j$ and $\sigma$ in $\sym_r$ let $\sigma_{ij}$ be the associated permutation in $\sym_2$.
Given $(\mu,\sigma) \in \CG(r)$ then an
element $w = (w_0,\dots w_n) \in \BE(r)_n$
$w \in \sur(r)_{(\mu,\sigma)}$ if for each $i<j$, the cardinality of $\set{\ell \mid (w_{\ell})_{ij} \neq (w_{\ell+1})_{ij}}$ is either less than $\mu_{ij}$ or equal to it and $(w_0)_{ij} = \sigma_{ij}$.

In particular $w$ has complexity $d$ or less if for each $i<j$
%the cardinality of $\set{\ell=1,\dots,n-1 \mid w_{\ell-1,ij} \neq w_{\ell,ij}}$ is less than $d$, i.e., 
the non-degenerate dimension of $w_{ij}=((w_0)_{ij},\dots,(w_n)_{ij})$ in $\BE(2)$ is $d$ or less for all $i<j$.
We notice that the action of $\sym_r$ on $\BE(r)$ preserves the nested sequence
$$\BE_1(r) \subset \BE_2(r) \subset \dots$$
For a proof that this is a $\CG(r)$-filtration we refer to Example 2.8 in \cite{berger1997confspacemodel}.
Please consult \cite{smith1989filtration,kashiwabara1993confcomplex,berger1997confspacemodel} for more details.
 
Applying the functor of singular chains to the zig-zag of \cref{p:zig-zag conf} produces a zig-zag of equivariant quasi-isomorphisms of $\UM$-coalgebras connecting $\schains\bars{\BE_d(r)}$ and $\schains\conf{r}{d}$.
Using the unit of the Quillen equivalence extends this zig-zag to one relating $\chains\BE_d(r)$ and $\schains\conf{r}{d}$, which can be combined with the zig-zag constructed in the previous subsection.
As announced in the introduction, this construction relates the chains on the multisimplicial model of configuration space and those the simplicial model via an explicit zig-zag of equivariant quasi-isomorphisms of $E_\infty$-coalgebras.

%% file: sec/tc.tex
% !TEX root = ../msimp.tex

\subsection{Table completion}\label{ss:table completion}

It is desirable to have a direct $\sym_r$-equivariant quasi-isomorphism between these algebraic models.
Two filtration preserving quasi-isomorphisms were constructed by Berger--Fresse
\[
\TR \colon \chains\BE(r) \to \chains\sur(r)
\quad \text{and} \quad
\TC \colon \chains\sur(r) \to \chains \BE(r).
\]
The first one, introduced in \cite[1$\cdot$3]{berger2004combinatorial}, is not a coalgebra map, as the reader familiar with its definition can easily verify.
We will focus on the second one which was introduced in \cite{berger2002prismatic} and termed \textit{table completion}.
We will construct a factorization
\[
\TC \colon \chains\sur(r) \xra{\EZ} \chains\sur(r)^\diag \xra{\chains(\tc)} \chains \BE(r),
\]
where the second map is induced from a simplicial map defined below.
This factorization proves that $\TC$ is a coalgebra map since both factors are.
We warn the reader that since $\EZ$ does not respect the $E_\infty$-coalgebra structure,
neither does $\TC$.
For example, we have
\[
\Delta_1 \circ \TC\big(12312) \neq \TC^{\otimes 2} \circ \, \Delta_1 \big(12312\big).
\]

Let us give a more explicit description of the simplicial set $\sur(r)^\diag$.
Its $n$-simplices are represented by maps
\[
f \colon \set{1,\dots,rm+r} \to \set{1,\dots,r},
\]
satisfying that the cardinality of $f^{-1}(\ell)$ is $m+1$ for each $\ell \in \set{1,\dots,r}$.
We represented this simplex by the sequence $f(1) \dotsb f(rm+r)$.
The $i^\th$ face and degeneracy maps act on it respectively by removing or doubling the $i^\th$ occurrence of each $\ell \in \set{1,\dots,r}$ in $f(1) \dotsb f(rm+r)$.

The restriction to the diagonal defines a $\CG(r)$-filtration of $\sur(r)^\diag$ and a nested sequence
\[
\sur_1(r)^\diag \subset \sur_2(r)^\diag \subset \dotsb
\]
on $\sur(r)^\diag$ that is preserved by the action of $\sym_r$ on $\sur(r)^\diag$.
In terms of cellular $\CG(r)$-decompositions, given $(\mu,\sigma) \in \CG(r)$ then $f \in \sur(r)^\diag_{(\mu,\sigma)}$ if for each $i<j$, either $i$ and $j$ alternate strictly less than $\mu_{ij}$ times in the sequence $f_{ij}$, or they do so exactly $\mu_{ij}$ times and the ordering formed by the first occurrences of $i$ and $j$ in $f_{ij}$ agrees with $\sigma_{ij}$.

Since the complexity of an element is unchanged by degeneracy maps, it can easily be seen that $\EZ \colon \chains\sur(r) \to \chains\sur(r)^\diag$ preserves $\CG(r)$-filtrations.

Let us now define the simplicial map $\tc$.
For $f$ as above, let
\[
\tc(f) = (\sigma_0,\dots,\sigma_m)
\]
with $\sigma_j$ represented by the subsequence of $f$ containing the $(j+1)^{\mathrm{st}}$ occurrence of each $\ell \in \{1,\dots,r\}$.
For example, we have
\[
\tc(122333112) = (123,231,312).
\]

\begin{theorem}
	The simplicial map $\tc \colon \sur(r)^\diag \to \BE(r)$ satisfies
	\[
	\TC = \chains(\tc) \circ \EZ
	\]
	and induces a weak equivalence
	\[
	\tc_d \colon \sur_d(r)^\diag \to \BE_d(r)
	\]
	for every $r,d \in \N$.
\end{theorem}

\begin{proof}
	It is clear that $\tc$ is a simplicial map, and verifying its relationship with Berger--Fresse's chain map is straightforward using the prism interpretation of a surjection as in \cite[\S2.1]{berger2002prismatic} and explicit combinatorial formula for $TC$ in \cite[\S3.1]{berger2002prismatic}.

	To check that $\tc$ preserves $\CG(r)$-filtrations let us first notice that $\tc(f)_{ij} = \tc(f_{ij})$, so without loss of generality we can assume $r=2$.
	In this case it is clear that non-degenerate simplices are sent to non-degenerate simplices (of the same dimension), so for these the complexity is preserved.
	We conclude the same for degenerate simplices using that $\tc$ is a simplicial map and that degeneracy maps leave complexity unchanged.

	For each $(\mu,\sigma) \in \CG_d(r)$, $\sur^\diag(r)_{(\mu,\sigma)}$ is mapped to $\BE(r)_{(\mu,\sigma)}$ since the order in which $i$ and $j$ first appear in $f_{ij}$ determines the first simplex in $\tc(f)_{ij}$.
	The final claim then follows the naturality of \eqref{e:cellular poset decomposition}.
\end{proof}

%% file: sec/count.tex
% !TEX root = ../msimp.tex

\subsection{Counting generators}

We would like to stress that the number of non-degenerate multisimplices in $\sur_d(r)$ is much smaller than the number of non-degenerate simplices in $\BE_d(r)$.
For example,
\begin{align*}
	& P_\chi^{2,4}(x) = 24(1+6x+10x^2+5x^3) \\
	& P_\cE^{2,4}(x) = 24(1+23x+104x^2+196x^3+184x^4+86x^5+16x^6)
\end{align*}
and
\begin{align*}
	& P_\chi^{3,3}(x) = 6(1+3x+7x^2+9x^3+6x^4+x^5) \\
	& P_\cE^{3,3}(x) = 6(1+5x+25x^2+60x^3+70x^4+38x^5+8x^6 )
\end{align*}
where
\begin{align*}
	P_\chi^{d,r}(x) &= \sum_n \,\rank(\chi_d(r)_n) \cdot x^n, \\
	P_\cE^{d,r}(x) &= \sum_n \,\rank(\cE_d(r)_n) \cdot x^n.
\end{align*}
This makes the multisimplicial approach substantially more efficient than the simplicial
when performing computations.
A calculation of obstruction to formality similar to that in
\cite{salvatore2020planarnonformality} took a full day with the simplicial model, and few seconds with the multisimplicial model.

%% file: aux/bibliography.bib
@article {brunfiedorowiczvogt2007,
    AUTHOR = {Brun, Morten and Fiedorowicz, Zbigniew and Vogt, Rainer M.},
     TITLE = {On the multiplicative structure of topological {H}ochschild
              homology},
   JOURNAL = {Algebr. Geom. Topol.},
  FJOURNAL = {Algebraic \& Geometric Topology},
    VOLUME = {7},
      YEAR = {2007},
     PAGES = {1633--1650},
      ISSN = {1472-2747},
   MRCLASS = {55P43 (18D50 19D55 55P48)},
  MRNUMBER = {2366174},
MRREVIEWER = {Vigleik Angeltveit},
       DOI = {10.2140/agt.2007.7.1633},
       URL = {https://doi.org/10.2140/agt.2007.7.1633},
}

@inproceedings {quillen2010higheralgktheory,
    AUTHOR = {Quillen, Daniel},
     TITLE = {Higher algebraic {$K$}-theory. {I}},
 BOOKTITLE = {Algebraic {$K$}-theory, {I}: {H}igher {$K$}-theories ({P}roc.
              {C}onf., {B}attelle {M}emorial {I}nst., {S}eattle, {W}ash.,
              1972)},
    SERIES = {Lecture Notes in Math., Vol. 341},
     PAGES = {85--147},
 PUBLISHER = {Springer, Berlin},
      YEAR = {1973},
   MRCLASS = {18F25},
  MRNUMBER = {0338129},
MRREVIEWER = {Stephen M. Gersten},
       DOI = {10.1007/BFb0067053},
       URL = {https://link.springer.com/chapter/10.1007/BFb0067053},
}

@article {smith1989filtration,
	AUTHOR = {Smith, Jeffrey Henderson},
	TITLE = {Simplicial group models for {$\Omega^nS^nX$}},
	JOURNAL = {Israel J. Math.},
	FJOURNAL = {Israel Journal of Mathematics},
	VOLUME = {66},
	YEAR = {1989},
	NUMBER = {1-3},
	PAGES = {330--350},
	ISSN = {0021-2172},
	MRCLASS = {55P47 (55U10)},
	MRNUMBER = {1017171},
	MRREVIEWER = {Ronald Brown},
	DOI = {10.1007/BF02765902},
	URL = {https://doi.org/10.1007/BF02765902},
}

@article {mcluresmith2004geomodel,
	AUTHOR = {McClure, James E. and Smith, Jeffrey H.},
	TITLE = {Cosimplicial objects and little {$n$}-cubes. {I}},
	JOURNAL = {Amer. J. Math.},
	FJOURNAL = {American Journal of Mathematics},
	VOLUME = {126},
	YEAR = {2004},
	NUMBER = {5},
	PAGES = {1109--1153},
	ISSN = {0002-9327},
	MRCLASS = {55P48 (18D50 55U10)},
	MRNUMBER = {2089084},
	MRREVIEWER = {Beno\^{\i}t Fresse},
	URL = {http://muse.jhu.edu/journals/american_journal_of_mathematics/v126/126.5mcclure.pdf},
}

@article {salvatore2009deligne,
	AUTHOR = {Salvatore, Paolo},
	TITLE = {The topological cyclic {D}eligne conjecture},
	JOURNAL = {Algebr. Geom. Topol.},
	FJOURNAL = {Algebraic \& Geometric Topology},
	VOLUME = {9},
	YEAR = {2009},
	NUMBER = {1},
	PAGES = {237--264},
	ISSN = {1472-2747},
	MRCLASS = {55P48 (18D50 55U10)},
	MRNUMBER = {2482075},
	MRREVIEWER = {Beno\^{\i}t Fresse},
	DOI = {10.2140/agt.2009.9.237},
	URL = {https://doi.org/10.2140/agt.2009.9.237},
}

@incollection {moerdijk1989bisimplicialsets,
	AUTHOR = {Moerdijk, Ieke},
	TITLE = {Bisimplicial sets and the group-completion theorem},
	BOOKTITLE = {Algebraic {$K$}-theory: connections with geometry and topology
	({L}ake {L}ouise, {AB}, 1987)},
	SERIES = {NATO Adv. Sci. Inst. Ser. C: Math. Phys. Sci.},
	VOLUME = {279},
	PAGES = {225--240},
	PUBLISHER = {Kluwer Acad. Publ., Dordrecht},
	YEAR = {1989},
	MRCLASS = {18G30 (55U10)},
	MRNUMBER = {1045852},
	URL = {https://doi.org/10.1007/978-94-009-2399-7_10},
	DOI = {10.1007/978-94-009-2399-7_10},
}

@article {salvatore2020planarnonformality,
	AUTHOR = {Salvatore, Paolo},
	TITLE = {Non-formality of planar configuration spaces in characteristic	2},
	JOURNAL = {Int. Math. Res. Not. IMRN},
	FJOURNAL = {International Mathematics Research Notices. IMRN},
	YEAR = {2020},
	NUMBER = {10},
	PAGES = {3100--3129},
	ISSN = {1073-7928},
	MRCLASS = {55R80 (55P62)},
	MRNUMBER = {4098635},
	MRREVIEWER = {Sadok Kallel},
	DOI = {10.1093/imrn/rny091},
	URL = {https://doi.org/10.1093/imrn/rny091},
}

@incollection {kashiwabara1993confcomplex,
	AUTHOR = {Kashiwabara, Takuji},
	TITLE = {On the homotopy type of configuration complexes},
	BOOKTITLE = {Algebraic topology ({O}axtepec, 1991)},
	SERIES = {Contemp. Math.},
	VOLUME = {146},
	PAGES = {159--170},
	PUBLISHER = {Amer. Math. Soc., Providence, RI},
	YEAR = {1993},
	MRCLASS = {55P35 (55S12)},
	MRNUMBER = {1224913},
	MRREVIEWER = {Hal Sadofsky},
	DOI = {10.1090/conm/146/01221},
	URL = {https://doi.org/10.1090/conm/146/01221},
}

@incollection {berger1997confspacemodel,
	AUTHOR = {Berger, Clemens},
	TITLE = {Combinatorial models for real configuration spaces and
	{$E_n$}-operads},
	BOOKTITLE = {Operads: {P}roceedings of {R}enaissance {C}onferences
	({H}artford, {CT}/{L}uminy, 1995)},
	SERIES = {Contemp. Math.},
	VOLUME = {202},
	PAGES = {37--52},
	PUBLISHER = {Amer. Math. Soc., Providence, RI},
	YEAR = {1997},
	MRCLASS = {18D35 (06A07 18B35 20B30 55P35 55U10)},
	MRNUMBER = {1436916},
	MRREVIEWER = {Peter J. Eccles},
	DOI = {10.1090/conm/202/02582},
	URL = {https://doi.org/10.1090/conm/202/02582},
}

@article {barrat1974operad,
	AUTHOR = {Barratt, M. G. and Eccles, Peter J.},
	TITLE = {{$\Gamma \sp{+}$}-structures. {I}. {A} free group functor for
	stable homotopy theory},
	JOURNAL = {Topology},
	FJOURNAL = {Topology. An International Journal of Mathematics},
	VOLUME = {13},
	YEAR = {1974},
	PAGES = {25--45},
	ISSN = {0040-9383},
	MRCLASS = {55D35 (55E10)},
	MRNUMBER = {348737},
	MRREVIEWER = {J. P. May},
	DOI = {10.1016/0040-9383(74)90036-6},
	URL = {https://doi.org/10.1016/0040-9383(74)90036-6},
}

@article {berger2002prismatic,
	AUTHOR = {Berger, Clemens and Fresse, Benoit},
	TITLE = {Une d\'{e}composition prismatique de l'op\'{e}rade de
	{B}arratt-{E}ccles},
	JOURNAL = {C. R. Math. Acad. Sci. Paris},
	FJOURNAL = {Comptes Rendus Math\'{e}matique. Acad\'{e}mie des Sciences. Paris},
	VOLUME = {335},
	YEAR = {2002},
	NUMBER = {4},
	PAGES = {365--370},
	ISSN = {1631-073X},
	MRCLASS = {55P48 (18D50)},
	MRNUMBER = {1931518},
	MRREVIEWER = {J. M. Boardman},
	DOI = {10.1016/S1631-073X(02)02489-5},
	URL = {https://doi.org/10.1016/S1631-073X(02)02489-5},
}

@article {maltsiniotis2005grothendieck,
	AUTHOR = {Maltsiniotis, Georges},
	TITLE = {La th\'{e}orie de l'homotopie de {G}rothendieck},
	JOURNAL = {Ast\'{e}risque},
	FJOURNAL = {Ast\'{e}risque},
	NUMBER = {301},
	YEAR = {2005},
	PAGES = {vi+140},
	ISSN = {0303-1179},
	MRCLASS = {18G55 (14F35 18G10 55U10 55U35)},
	MRNUMBER = {2200690},
	MRREVIEWER = {J\v{i}r\'{\i} Rosick\'{y}},
	URL = {http://www.numdam.org/item/AST_2005__301__R1_0.pdf},
}


%% file: aux/usualpapers.bib
@article{medina2022cube_einfty,
	author = {{Kaufmann}, Ralph M. and {Medina-Mardones}, Anibal M.},
	title = "{A combinatorial ${E}_\infty$-algebra structure on cubical cochains and the Cartan--Serre map}",
	journal = {arXiv e-prints},
	year = {2022},
	archivePrefix = {arXiv},
	eprint = {2107.00669},
	doi = {10.48550/arXiv.2107.00669},
	note = {To appear in Cahiers Topologie G\'{e}om. Diff\'{e}rentielle Cat\'{e}g.},
}

@ARTICLE{medina2021comch,
	author = {Anibal M. {Medina-Mardones}},
	title = {{A computer algebra system for the study of commutativity up to coherent homotopies}},
	volume = {14},
	journal = {Advanced Studies: Euro-Tbilisi Mathematical Journal},
	number = {4},
	publisher = {Tbilisi Centre for Mathematical Sciences},
	pages = {147 -- 157},
	year = {2021},
	url = {https://projecteuclid.org/journals/advanced-studies-euro-tbilisi-mathematical-journal/volume-14/issue-4/A-computer-algebra-system-for-the-study-of-commutativity-up/10.3251/asetmj/1932200819.full}
}

@ARTICLE{medina2021prop2,
	AUTHOR = {{Medina-Mardones}, Anibal M.},
	TITLE = {A finitely presented {$E_\infty$}-prop {II}: cellular context},
	JOURNAL = {High. Struct.},
	FJOURNAL = {Higher Structures},
	VOLUME = {5},
	YEAR = {2021},
	NUMBER = {1},
	PAGES = {69--186},
	URL = {https://higher-structures.math.cas.cz/api/files/issues/Vol5Iss1/Medina-Mardones-2}
}

@article {medina2021may_st,
	AUTHOR = {Kaufmann, Ralph M. and Medina-Mardones, Anibal M.},
	TITLE = {Cochain level {M}ay-{S}teenrod operations},
	JOURNAL = {Forum Math.},
	FJOURNAL = {Forum Mathematicum},
	VOLUME = {33},
	YEAR = {2021},
	NUMBER = {6},
	PAGES = {1507--1526},
	ISSN = {0933-7741},
	MRCLASS = {55S05 (55S10 55S12 55S15 55U05 55U15)},
	MRNUMBER = {4333989},
	DOI = {10.1515/forum-2020-0296},
	URL = {https://doi.org/10.1515/forum-2020-0296},
}

@article {medina2021adem,
	AUTHOR = {Brumfiel, Greg and {Medina-Mardones}, Anibal and Morgan, John},
	TITLE = {A cochain level proof of {A}dem relations in the mod 2 {S}teenrod algebra},
	JOURNAL = {J. Homotopy Relat. Struct.},
	FJOURNAL = {Journal of Homotopy and Related Structures},
	VOLUME = {16},
	YEAR = {2021},
	NUMBER = {4},
	PAGES = {517--562},
	ISSN = {2193-8407},
	MRCLASS = {55S10},
	MRNUMBER = {4343073},
	DOI = {10.1007/s40062-021-00287-3},
	URL = {https://doi.org/10.1007/s40062-021-00287-3},
}

@ARTICLE{medina2020prop1,
	AUTHOR = {{Medina-Mardones}, Anibal M.},
	TITLE = {A finitely presented {$E_\infty$}-prop {I}: algebraic context},
	JOURNAL = {High. Struct.},
	FJOURNAL = {Higher Structures},
	VOLUME = {4},
	YEAR = {2020},
	NUMBER = {2},
	PAGES = {1--21},
	MRCLASS = {55P48 (18M85 18N50)},
	MRNUMBER = {4133162},
	url = {https://journals.mq.edu.au/api/files/issues/Vol4Iss2/Medina-Mardones}
}

@ARTICLE{medina2020globular,
	AUTHOR = {{Medina-Mardones}, Anibal M.},
	TITLE = {An algebraic representation of globular sets},
	JOURNAL = {Homology Homotopy Appl.},
	FJOURNAL = {Homology, Homotopy and Applications},
	VOLUME = {22},
	YEAR = {2020},
	NUMBER = {2},
	PAGES = {135--150},
	ISSN = {1532-0073},
	MRCLASS = {18N99 (18F99 55S05)},
	MRNUMBER = {4093174},
	DOI = {10.4310/hha.2020.v22.n2.a8},
	URL = {https://doi.org/10.4310/hha.2020.v22.n2.a8},
}

@ARTICLE{medina2020cartan,
	AUTHOR = {{Medina-Mardones}, Anibal M.},
	TITLE = {An effective proof of the {C}artan formula: the even prime},
	JOURNAL = {J. Pure Appl. Algebra},
	FJOURNAL = {Journal of Pure and Applied Algebra},
	VOLUME = {224},
	YEAR = {2020},
	NUMBER = {12},
	PAGES = {106444, 18},
	ISSN = {0022-4049},
	MRCLASS = {55S10 (55S05 55S12)},
	MRNUMBER = {4102178},
	DOI = {10.1016/j.jpaa.2020.106444},
	URL = {https://doi.org/10.1016/j.jpaa.2020.106444},
}

@ARTICLE{medina2021cobar,
	author = {{Medina-Mardones}, Anibal M. and {Rivera}, Manuel},
	title = "{The cobar construction as an $E_{\infty}$-bialgebra model of the based loop space}",
	journal = {arXiv e-prints},
	year = {2021},
	archivePrefix = {arXiv},
	eprint = {2108.02790},
	note = {Submitted}
}

@ARTICLE{medina2022axiomatic,
	author = {{Medina-Mardones}, Anibal M.},
	title = "{An axiomatic characterization of {S}teenrod's cup-$i$ products}",
	journal = {arXiv e-prints},
	year = {2022},
	archivePrefix = {arXiv},
	eprint = {1810.06505},
	note = {Submitted},
}

@article {mcclure2003multivariable,
	AUTHOR = {McClure, James E. and Smith, Jeffrey H.},
	TITLE = {Multivariable cochain operations and little {$n$}-cubes},
	JOURNAL = {J. Amer. Math. Soc.},
	FJOURNAL = {Journal of the American Mathematical Society},
	VOLUME = {16},
	YEAR = {2003},
	NUMBER = {3},
	PAGES = {681--704},
	ISSN = {0894-0347},
	MRCLASS = {55P48 (18D50)},
	MRNUMBER = {1969208},
	MRREVIEWER = {Benoit Fresse},
	DOI = {10.1090/S0894-0347-03-00419-3},
	URL = {https://doi.org/10.1090/S0894-0347-03-00419-3},
}

@article {berger2004combinatorial,
	AUTHOR = {Berger, Clemens and Fresse, Benoit},
	TITLE = {Combinatorial operad actions on cochains},
	JOURNAL = {Math. Proc. Cambridge Philos. Soc.},
	FJOURNAL = {Mathematical Proceedings of the Cambridge Philosophical
	Society},
	VOLUME = {137},
	YEAR = {2004},
	NUMBER = {1},
	PAGES = {135--174},
	ISSN = {0305-0041},
	MRCLASS = {18D50 (16E45 55P48)},
	MRNUMBER = {2075046},
	MRREVIEWER = {David Chataur},
	DOI = {10.1017/S0305004103007138},
	URL = {https://doi.org/10.1017/S0305004103007138},
}

@article {mandell2006homotopy_type,
	AUTHOR = {Mandell, Michael A.},
	TITLE = {Cochains and homotopy type},
	JOURNAL = {Publ. Math. Inst. Hautes \'{E}tudes Sci.},
	FJOURNAL = {Publications Math\'{e}matiques. Institut de Hautes \'{E}tudes
	Scientifiques},
	NUMBER = {103},
	YEAR = {2006},
	PAGES = {213--246},
	ISSN = {0073-8301},
	MRCLASS = {55P15 (55N10 55P60)},
	MRNUMBER = {2233853},
	MRREVIEWER = {Donald Yau},
	DOI = {10.1007/s10240-006-0037-6},
	URL = {https://doi.org/10.1007/s10240-006-0037-6},
}

@article {steenrod1947products,
	AUTHOR = {Steenrod, N. E.},
	TITLE = {Products of cocycles and extensions of mappings},
	JOURNAL = {Ann. of Math. (2)},
	FJOURNAL = {Annals of Mathematics. Second Series},
	VOLUME = {48},
	YEAR = {1947},
	PAGES = {290--320},
	ISSN = {0003-486X},
	MRCLASS = {56.0X},
	MRNUMBER = {22071},
	MRREVIEWER = {B. Eckmann},
	DOI = {10.2307/1969172},
	URL = {https://doi.org/10.2307/1969172},
}

@article {steenrod1952reduced,
	AUTHOR = {Steenrod, N. E.},
	TITLE = {Reduced powers of cohomology classes},
	JOURNAL = {Ann. of Math. (2)},
	FJOURNAL = {Annals of Mathematics. Second Series},
	VOLUME = {56},
	YEAR = {1952},
	PAGES = {47--67},
	ISSN = {0003-486X},
	MRCLASS = {56.0X},
	MRNUMBER = {48026},
	MRREVIEWER = {H. Cartan},
	DOI = {10.2307/1969766},
	URL = {https://doi.org/10.2307/1969766},
}

@article {steenrod1953cyclic,
	AUTHOR = {Steenrod, N. E.},
	TITLE = {Cyclic reduced powers of cohomology classes},
	JOURNAL = {Proc. Nat. Acad. Sci. U.S.A.},
	FJOURNAL = {Proceedings of the National Academy of Sciences of the United
	States of America},
	VOLUME = {39},
	YEAR = {1953},
	PAGES = {217--223},
	ISSN = {0027-8424},
	MRCLASS = {56.0X},
	MRNUMBER = {54965},
	MRREVIEWER = {H. Cartan},
	DOI = {10.1073/pnas.39.3.217},
	URL = {https://doi.org/10.1073/pnas.39.3.217},
}

@article {gugenheim1957supercomplexes,
	AUTHOR = {Gugenheim, V. K. A. M.},
	TITLE = {On supercomplexes},
	JOURNAL = {Trans. Amer. Math. Soc.},
	FJOURNAL = {Transactions of the American Mathematical Society},
	VOLUME = {85},
	YEAR = {1957},
	PAGES = {35--51},
	ISSN = {0002-9947},
	MRCLASS = {55.0X},
	MRNUMBER = {86299},
	MRREVIEWER = {J.-P. Meyer},
	DOI = {10.2307/1992960},
	URL = {https://doi.org/10.2307/1992960},
}
